# RANDOM FIELDS OF MULTIVARIATE TEST STATISTICS, WITH APPLICATIONS TO SHAPE ANALYSIS

By J. E. Taylor[1] and K. J. Worsley[2]

*Stanford University, Université de Montréal and McGill University*

Our data are random fields of multivariate Gaussian observations, and we fit a multivariate linear model with common design matrix at each point. We are interested in detecting those points where some of the coefficients are nonzero using classical multivariate statistics evaluated at each point. The problem is to find the $P$-value of the maximum of such a random field of test statistics. We approximate this by the expected Euler characteristic of the excursion set. Our main result is a very simple method for calculating this, which not only gives us the previous result of Cao and Worsley [*Ann. Statist.* **27** (1999) 925–942] for Hotelling's $T^2$, but also random fields of Roy's maximum root, maximum canonical correlations [*Ann. Appl. Probab.* **9** (1999) 1021–1057], multilinear forms [*Ann. Statist.* **29** (2001) 328–371], $\bar{\chi}^2$ [*Statist. Probab. Lett* **32** (1997) 367–376, *Ann. Statist.* **25** (1997) 2368–2387] and $\chi^2$ scale space [*Adv. in Appl. Probab.* **33** (2001) 773–793]. The trick involves approaching the problem from the point of view of Roy's union-intersection principle. The results are applied to a problem in shape analysis where we look for brain damage due to nonmissile trauma.

**1. Introduction.** Our motivation comes from a study by Tomaiuolo, Worsley, Lerch, Di Paulo, Carlesimo, Bonanni, Caltagirone and Paus [25] on the anatomy of a group of 17 nonmissile brain trauma patients measured by magnetic resonance imaging (MRI). The aim is to detect regions of brain damage (shape change) by comparing anatomy at each point in $N = 3$ dimensional space to that of a group of 19 age and sex matched controls. Each brain was first linearly transformed into a common stereotactical reference space. Then the method of Collins, Holmes, Peters and Evans [6] was used

Received July 2005; revised March 2007.
[1]Supported in part by NSF Grant DMS-04-05970.
[2]Supported in part by NSERC and FQRNT.
*AMS 2000 subject classifications.* 52A22, 62H11.
*Key words and phrases.* Integral geometry, differential topology, Euler characteristic, Morse theory, Roy's maximum root, Hotelling's $T^2$, canonical correlation, scale space.

to find the nonlinear vector deformations $Y_i(s) \in \Re^d$ ($d=3$ here) required to transform the MRI image of subject $i$ to a common atlas standard at each point $s$ inside a *search region* $S \subset \Re^N$, usually the whole brain—see Figure 1. This sort of anatomical data is good at detecting changes in the boundary of brain structures, with the added benefit of estimating the direction in which the change took place.

To do this, Cao and Worsley [4] set up a linear model for subject $i$, $i = 1, \ldots, n$:

$$Y_i(s)' = x_i'\beta(s) + Z_i(s)'\Sigma(s)^{1/2}, \tag{1}$$

where $x_i$ is a $p$-vector of known regressors, and $\beta(s)$ is an unknown $p \times d$ coefficient matrix. The error $Z_i(s)$ is a $d$-vector of independent zero mean, unit variance Gaussian random field components with the same spatial correlation structure, and $\text{Var}(Y_i(s)) = \Sigma(s)$ is an unknown $d \times d$ matrix. We can now detect how the regressors are related to shape at point $s$ by testing contrasts in the rows of $\beta(s)$. Classical multivariate test statistics evaluated at each point $s$ then form a random field $T(s)$.

We suspect that the changes are confined to a small number of localized regions and the rest of the brain is unchanged. It is not hard to show that if the spatial pattern of the change matches the spatial correlation function of the errors, and if $T(s)$ is the likelihood ratio at a single point $s$, then the spatial maximum of $T(s)$ is the likelihood ratio test under unknown signal location, provided $\Sigma(s)$ is known (Siegmund and Worsley [14]). In essence, this is just the Matched Filter Theorem from signal processing. Thresholding $T(s)$ at that threshold which controls the $P$-value of its maximum should then be powerful at detecting changes, while controlling the false positive rate outside the changed region to something slightly smaller than the nominal $P$-value. Our main problem is therefore to find the $P$-value of the maximum of random fields of multivariate test statistics.

The outline of the paper is as follows. An approximate $P$-value of the maximum of a random field is given in Section 2, and in Section 3 we apply this to random fields of multivariate test statistics $T(s)$ such as Hotelling's $T^2$, Roy's maximum root and maximum canonical correlation. The same methods can also be applied to random fields of $\bar{\chi}^2$ (Lin and Lindsay [12] and Takemura and Kuriki [18]) statistics, but this will be the subject of a forthcoming paper (Taylor and Worsley [23]). In Section 4 we apply these results to the nonmissile brain trauma data above. A further application to $\chi^2$ scale space random fields is given in Appendix A.5.

## 2. *P*-value of the maximum of a random field.

2.1. *The expected Euler characteristic.* A random field $T(s)$, $s \in S \subset \mathbb{R}^N$, is *isotropic* if it has the same distribution as $T(a + Bs)$, where $a \in \mathbb{R}^N$ is any translation and $B$ is any orthnormal (rotation) matrix. A very accurate approximation to the $P$-value of the maximum of such a random field, at high thresholds $t$, is the expected Euler characteristic (EC) $\varphi$ of the excursion set

$$(2) \qquad \mathbb{P}\Big(\max_{s \in S} T(s) \geq t\Big) \approx \mathbb{E}(\varphi\{s \in S : T(s) \geq t\}) = \sum_{i=0}^{N} \mu_i(S) \rho_i(t),$$

where $\mu_i(S)$ is the $i$-dimensional *intrinsic volume* of $S$, and $\rho_i(t)$ is the $i$-dimensional *EC density* of the random field above $t$ (Adler [1, 2], Worsley [27], Taylor, Takemura and Adler [22] and Adler and Taylor [3]). The accuracy of the approximation (2) will be discussed in Section 2.2.

For $N = 3$, our main interest in applications, the EC of a set is

$$\varphi = \sharp\text{blobs} - \sharp\text{handles or tunnels} + \sharp\text{interior hollows}$$

(see Figure 2). It can be evaluated numerically for a subset of a rectilinear mesh by

$$(3) \qquad \varphi = \sharp\text{points} - \sharp\text{edges} + \sharp\text{faces} - \sharp\text{cubes},$$

where, for example, a cube is a set of 8 adjacent mesh points, differing by one mesh step along each axis, and all inside the set [Figure 2(a)]. This method was used to calculate the EC of the excursion sets in Figure 2(b)–(e).

Intrinsic volumes are defined in Appendix A.1, and for $N = 3$, they are

$$(4) \qquad \begin{aligned} \mu_0(S) &= \varphi(S), \\ \mu_1(S) &= 2 \text{ caliper diameter}(S), \\ \mu_2(S) &= \tfrac{1}{2}\text{surface area}(S), \\ \mu_3(S) &= \text{volume}(S). \end{aligned}$$

For convex $S$, the caliper diameter is the distance between two parallel planes tangent to $S$, averaged over all rotations.

EC density is defined in Section 2.3. As an example, suppose $T = T(s)$ is a *unit Gaussian random field* (UGRF), defined as a Gaussian random field with $\mathbb{E}(T) = 0$, $\text{Var}(T) = 1$ and $\text{Var}(\dot{T}) = I_{N \times N}$, the $N \times N$ identity matrix, where dot denotes (spatial) derivative with respect to $s$. Then the first four

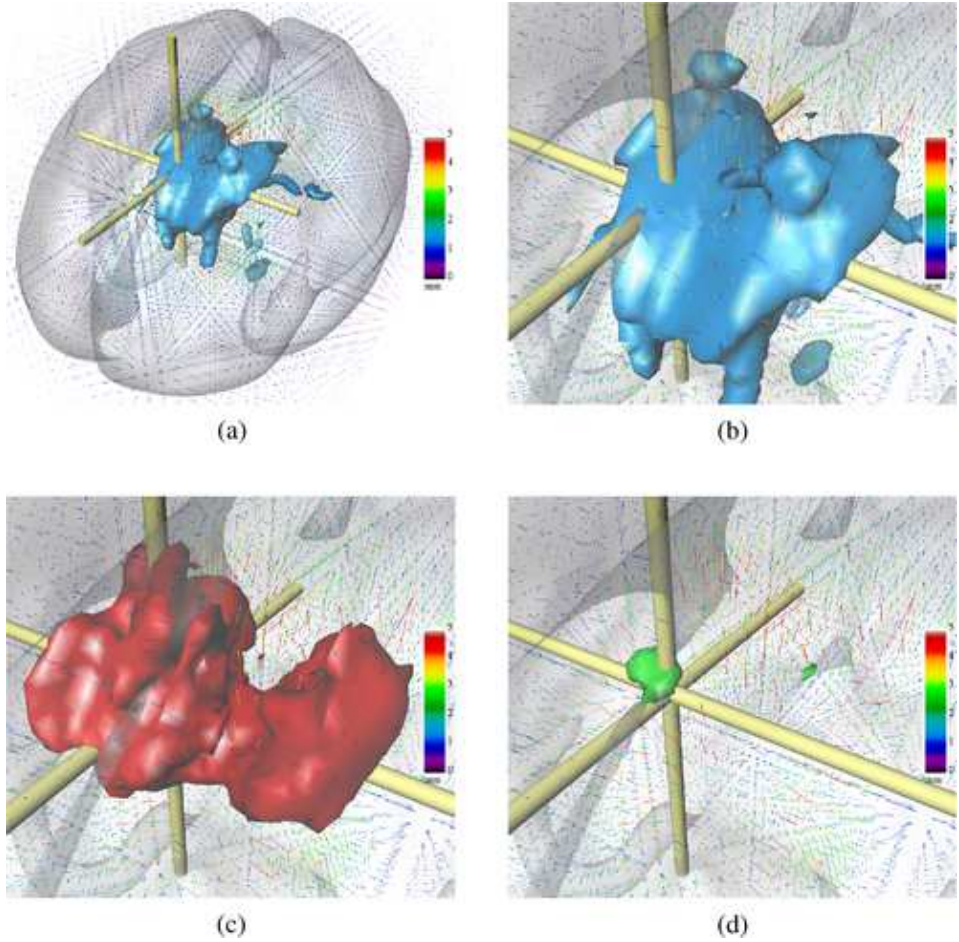

FIG. 1. *Shape analysis of nonmissile brain trauma data.* (a) *Trauma minus control average deformations (arrows and color bar), sampled every 6mm inside the brain, with $T(s) =$ Hotelling's $T^2$ field for significant group differences (threshold $t = 54.0$, $P = 0.05$). The reference point of maximum Hotelling's $T^2$ is marked by the intersection of the three axes.* (b) *Closeup of* (a) *showing that the damage is an outward movement of the anatomy, either due to swelling of the ventricles or atrophy of the surrounding white matter.* (c) *Regions of effective anatomical connectivity with the reference point, assessed by $T(s) =$ maximum canonical correlation field (threshold $t = 0.746$, $P = 0.05$). The reference point is connected with its neighbors (due to smoothness) and with contralateral regions (due to symmetry).* (d) *Regions where the connectivity is different between trauma and control groups, assessed by $T(s) =$ Roy's maximum root field (threshold $t = 30.3$, $P = 0.05$). The small region in the contralateral hemisphere (right) is more correlated in the trauma group than the control group.*

EC densities of $T$ are
$$\rho_0(t) = \mathbb{P}(T \geq t),$$

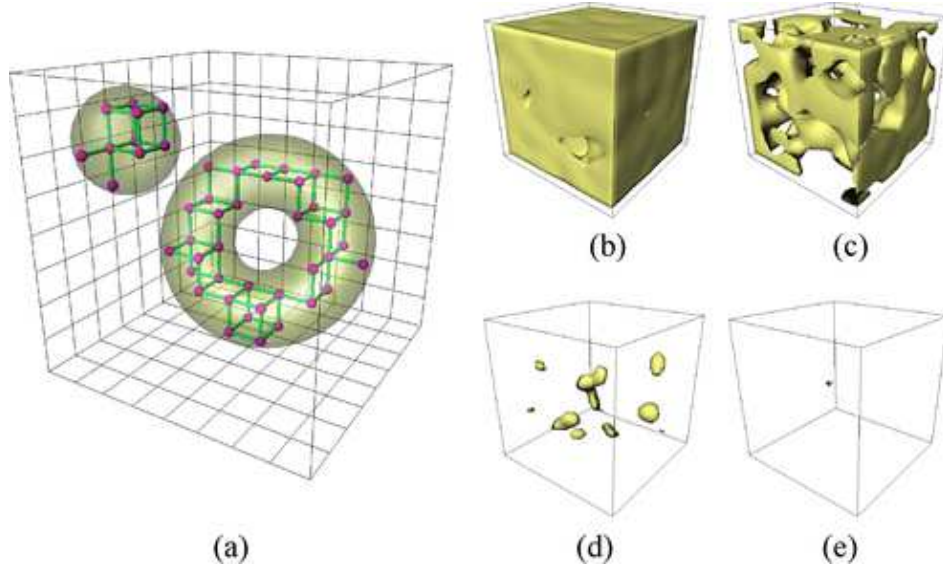

FIG. 2. (a) The Euler characteristic of a 3D set is $\varphi = \sharp blobs - \sharp handles\ or\ tunnels + \sharp hollows = 2 - 1 + 0 = 1$, or, if we fill the set with a fine rectilinear mesh, then $\varphi = \sharp points - \sharp edges + \sharp faces - \sharp cubes = 47 - 70 + 26 - 2 = 1$. (b) The excursion set of a 3D isotropic Gaussian random field with mean zero and variance one above a threshold $t = -2$; the set contains isolated hollows that each contribute $+1$ to give $\varphi = 5$ with $\mathbb{E}(\varphi) = 6.7$ from (2); (c) at $t = 0$ the handles or tunnels dominate, each contributing $-1$ to give $\varphi = -28$ with $\mathbb{E}(\varphi) = -20.0$; (d) at $t = 2$ the handles and hollows disappear, leaving isolated blobs, each contributing $+1$ to give $\varphi = 11$ with $\mathbb{E}(\varphi) = 16.1$; (e) at $t = 3$ only one blob remains (containing the maximum value of 3.16) to give $\varphi = 1$ with $\mathbb{E}(\varphi) = 2.1$. At very high thresholds $\mathbb{E}(\varphi)$ is a good approximation to the P-value of the maximum.

$$
\begin{aligned}
\rho_1(t) &= (2\pi)^{-1} \exp(-t^2/2), \\
\rho_2(t) &= (2\pi)^{-3/2} t \exp(-t^2/2), \\
\rho_3(t) &= (2\pi)^{-2} (t^2 - 1) \exp(-t^2/2).
\end{aligned}
\tag{5}
$$

Note that any stationary Gaussian random field can be transformed to a UGRF by an appropriate linear transformation of its domain and range. In particular, if $\text{Var}(\dot{T}) = cI_{N \times N}$ for some scalar $c$, then $\rho_i(t)$ is multiplied by $c^{i/2}$.

2.2. *The accuracy of the approximation.* A heuristic explanation for why we use the expected EC as a P-value approximation is as follows. If the threshold $t$ is high, the handles and hollows of the excursion set tend to disappear, leaving a set of isolated blobs, each containing one local maximum, so that the EC then counts the number of connected components [Figure 2(d)]. At very high thresholds the excursion set is mostly empty with an EC of 0, or occasionally, when $\max_{s \in S} T(s) \geq t$, it will contain just one connected component with an EC of 1 [Figure 2(e)]. Thus, at these high thresholds the expected EC is a good approximation to the $P$-value of $\max_{s \in S} T(s)$. The beauty of the EC is that there is an *exact* expression for its expectation for *all* thresholds.

Moreover, the approximation (2) is astonishingly accurate when $S$ is either convex or has smooth boundary, and $T(s)$ is Gaussian, in which case

$$
\begin{aligned}
(6) \quad & \mathbb{E}(\varphi\{s \in S : T(s) \geq t\}) \\
& = \mathbb{P}(T \geq t) + (c_0 + c_1 t + \cdots + c_{D-1} t^{D-1}) \exp(-t^2/2)
\end{aligned}
$$

for some constants $c_0, \ldots, c_{D-1}$ [see (5) and (9) below]. Since $\mathbb{P}(T \geq t) = O(1/t) \exp(-t^2/2)$, it might be thought that the error in the $P$-value approximation (2) is simply the next term down, that is, $O(1/t^2) \exp(-t^2/2)$. In fact, the error is *exponentially* smaller than this:

$$
(7) \qquad \mathbb{P}\left(\max_{s \in S} T(s) \geq t\right) = \mathbb{E}(\varphi\{s \in S : T(s) \geq t\}) + O(\exp(-\alpha t^2/2))
$$

for some $\alpha > 1$ which is related to the curvature of the boundary of $S$ and $\mathrm{Var}(\ddot{T})$ (Taylor, Takemura and Adler [22]). This means that there are in effect no further terms in the expansion, and the expected EC captures essentially all of the polynomial terms in the $P$-value.

We will see later that many of the $P$-value approximations for non-Gaussian fields can be transformed into $P$-value approximations of Gaussian fields on larger parameter spaces. This means that these special non-Gaussian fields also have exponentially accurate $P$-value approximations as in (7). There will be some examples, however, with a $\chi^2$ field in the denominator, for which this trick will not work. For these examples, it is difficult to give quantitative bounds on the error. The general techniques in Taylor, Takemura and Adler [22] still apply to these non-Gaussian fields, though the final form of the bound is not explicitly known in these cases.

2.3. *The EC density.* A direct method for finding EC densities of any sufficiently smooth random field can be derived using Morse theory (Morse and Cairns [13]) for $i > 0$,

$$
(8) \qquad \rho_i(t) = \mathbb{E}(1_{\{T \geq t\}} \det(-\ddot{T}_i) \mid \dot{T}_i = 0) \mathbb{P}(\dot{T}_i = 0),
$$

where dot notation with subscript $i$ denotes differentiation with respect to the first $i$ components of $s$, and double dot with subscript $i$ denotes the matrix of second derivatives with respect to the first $i$ components of $s$ (Adler [1] and Worsley [27]). For $i = 0$, $\rho_0(t) = \mathbb{P}(T \geq t)$. If $T$ is a UGRF, Adler [1] uses this method to obtain its EC density

$$\rho_i^{\mathrm{G}}(t) = \left( \frac{-1}{\sqrt{2\pi}} \frac{\partial}{\partial t} \right)^i \mathbb{P}(T \geq t), \tag{9}$$

which leads directly to (5) and (6).

Our main interest in this paper is finding $P$-values for maxima of random fields of test statistics commonly encountered in multivariate analysis, with the ultimate aim of applying this to detecting points $s$ where coefficients in a multivariate linear model are nonzero. First of all we must define the random field $T$. For example, a $\chi^2$ random field with $d$ degrees of freedom is defined as

$$T(s) = Z(s)'Z(s), \tag{10}$$

where $Z(s)$ is a $d$-vector of i.i.d. UGRFs. In a similar way we can define other test statistic random fields, such as $t$ and $F$ statistic random fields, by applying the usual definition to component UGRFs.

2.4. *The Gaussian kinematic formula.* Finding workable expressions for EC densities is a tedious process. The result for UGRFs (9) appears to be simple, but its derivation takes most of Chapter 5 of Adler [1]. In fact, Adler did not at first recognize that it could be expressed as the derivative of the $P$-value at a point. For many years this was thought to be a coincidence, but a new result of Taylor [20] on Gaussian kinematic formulas shows that indeed this is not a coincidence in certain cases (the Gaussian is one, the $\sqrt{\chi^2}$ is another): the EC densities of functions of UGRFs fall out as the coefficients of a power series expansion of the probability content of a tube about the rejection region, in terms of the tube radius. For the Gaussian case, this tube is just another rejection region of the same shape, which ultimately leads to the special form of the Gaussian EC density (9).

The details of the Gaussian kinematic formula are as follows. Suppose $T(s) = f(Z(s))$ is a function of $Z(s)$, a $d$-vector of i.i.d. UGRFs as before, so that the rejection or critical region is $C = \{z : f(z) \geq t\} \subset \mathbb{R}^d$. Let $\mathrm{Tube}(C, \varepsilon) = \{x : \min_{z \in C} \|z - x\| \leq \varepsilon\}$ be the tube of radius $\varepsilon$ around $C$. Let $Z$ be a $d$-vector of i.i.d. $\mathrm{N}(0,1)$ random variables. Then the Gaussian kinematic formula is

$$\mathbb{P}(Z \in \mathrm{Tube}(C, \varepsilon)) = \sum_{i=0}^{\infty} \frac{\varepsilon^i}{i!} (2\pi)^{i/2} \rho_i(t). \tag{11}$$

For $T(s)$ a UGRF, $d = 1$, $f(z) = z$, and this leads directly to the Gaussian EC density (9). For an $F$ statistic random field with $(\eta, \nu)$ degrees of freedom,

$$f(z) = \frac{z_1'z_1/\eta}{z_2'z_2/\nu},$$

where $z = (z_1', z_2')'$, $z_1 \in \mathbb{R}^\eta$, $z_2 \in \mathbb{R}^\nu$ and $d = \eta + \nu$. The rejection region $C$ is a cone, and after a little elementary geometry, it can be seen that the first $d$ terms come only from the sides of the cone and are easy to calculate (Taylor and Worsley [23]):

$$(12) \quad \mathbb{P}(Z \in \text{Tube}(C, \varepsilon)) = \mathbb{P}(\sqrt{V} \geq \sqrt{Wt\eta/\nu} - \varepsilon\sqrt{1 + t\eta/\nu}) + O(\varepsilon^d),$$

where $V \sim \chi^2_\eta$ independently of $W \sim \chi^2_\nu$. Expanding the right-hand side of (12) in powers of $\varepsilon$ and equating it to that of the Gaussian kinematic formula (11) gives the $i$-dimensional EC densities $\rho_i^F(t)$ of the $F$ field for $i < d$. This agrees with an explicit expression found using the Morse theory result (8) that can be found in Worsley [26].

The higher-order terms in $O(\varepsilon^d)$ in (12) come from the apex of the cone. Luckily we do not have to calculate them, since we only need EC densities up to dimension $N$ and if $N \geq d$, then the $F$ field is not defined (Worsley [26]). The $F$ field is not defined because both numerator and denominator of $F$ can take the value 0 (i.e., $F = 0/0$) with positive probability. The reason is that the zero set of any of the component UGRFs is a smooth surface of dimension $N - 1$; the intersection of $d$ such surfaces is a set of dimension $N - d$ which is nonempty (with positive probability) if $N - d \geq 0$. On this set both numerator and denominator of $F$ take the value 0 with positive probability and $F$ is not defined.

2.5. *Non-Gaussian random fields.* Unfortunately the EC densities of non-Gaussian random fields are not always derivatives of the $P$-value as in (9). Instead we must go back to the Morse theory result (8) or, if the random field is a function of i.i.d. UGRFs (as is the case here), we can use the Gaussian kinematic formula (11). Both are equally complex for the types of random fields we are interested in.

For the Gaussian kinematic formula, we must evaluate the probability content of a tube about a very complex rejection region. For the Morse theory result, we must obtain expressions for the joint distribution of the random field and its first two derivatives, then evaluate the expectation (8) by careful manipulation of vector and matrix random variables. This has been done for $\chi^2$, $t$ and $F$ random fields (Worsley [26]), Hotelling's $T^2$ random fields (Cao and Worsley [4]) and correlation random fields (Cao and Worsley [5]). In each case, after a great deal of tedious algebra, an exact closed-form expression for the EC density is obtained.

The types of random fields we are interested in are generalizations of Hotelling's $T^2$, which is used to test for a single coefficient in a multivariate linear model. Our goal is to generalize random field theory to testing for multiple coefficients, the only missing piece. The obvious choice is the likelihood ratio statistic, Wilks' $\Lambda$, but so far this has resisted our attempts. Instead, by a very simple trick that builds on previous EC densities and avoids any further evaluation of the Morse theory result (8), we can easily obtain $P$-value approximations similar to (2) (but not EC densities) for Roy's maximum root, a common alternative to Wilks' $\Lambda$. Using the same method, we re-derive in Section 3.1 the EC density for Hotelling's $T^2$, shortening the original derivation from an entire *Annals* paper (Cao and Worsley [4]) down to just several lines. We use the same trick to get the $\chi^2$ scale space field from the Gaussian scale space field, again reducing most of an *Advances in Applied Probability* paper (Worsley [28]) down to just one line—see Appendix A.5.

2.6. *The union-intersection principle.* The trick is to approach the problem from the point of view of Roy's union-intersection principle. Take, for example, the $\chi^2$ random field with $d$ degrees of freedom (10). We can write this as the maximum of linear combinations of $d$ i.i.d. UGRFs $Z(s)$ as follows:

$$T(s) = \max_u \tilde{T}(s,u), \qquad \text{where } \tilde{T}(s,u) = (Z(s)'u)^2,$$

and $u \in \mathbb{R}^d$ is restricted to the unit sphere $U_d = \{u : u'u = 1\}$. In other words, we have written the $\chi^2$ field in terms of the square of a Gaussian field over a larger domain $S \times U_d$.

At first glance this seems to make things more difficult, but, in fact, it makes things easier. We can now apply the $P$-value approximation (2) to $\tilde{T}(s,u)$, replacing $S$ by $S \times U_d$. Using integral geometry, the intrinsic volume of $S \times U_d$ is a simple function of the intrinsic volumes of $S$ and $U_d$. The EC densities of $\tilde{T}(s,u)$ are easily obtained from the Gaussian EC densities (9). Putting them all together and pulling out the coefficient of $\mu_i(S)$ gives us the EC density for the $\chi^2$ field.

2.7. *Nonisotropic fields and Lipschitz–Killing curvature.* There is only one technical difficulty. Although $\tilde{T}(s,u)$ is isotropic in $s$ for fixed $u$, it is not isotropic in $u$ for fixed $s$ (since $u$ is on a sphere), nor jointly isotropic in $(s,u)$. To handle this situation, Taylor and Adler [21] and Taylor [20] have extended the theory of Adler [1] and Worsley [27] to nonisotropic fields on smooth manifolds with piece-wise smooth boundaries. They show that if the random field is a function of i.i.d. nonisotropic UGRFs, then it is only necessary to replace intrinsic volume $\mu_i(S)$ in the expected EC (2) by *Lipschitz–Killing curvature* (LKC) $\mathcal{L}_i(S)$:

$$(13) \qquad \mathbb{P}\left(\max_{s \in S} T(s) \geq t\right) \approx \mathbb{E}(\varphi\{s \in S : T(s) \geq t\}) = \sum_{i=0}^{N} \mathcal{L}_i(S)\rho_i(t).$$

This result holds for nearly all nonisotropic random fields that are $N(0,1)$ at each point and smooth, that is, with at least two almost sure derivatives and an additional mild regularity condition; see Taylor and Adler [21].

The LKC of $S$ is a measure of the intrinsic volume of $S$ in the Riemannian metric defined by the variogram. Specifically, we replace local Euclidean distance between points $s_1, s_2 \in S$ by the square root of the variogram

$$\text{Var}(Z_1(s_1) - Z_1(s_2))^{1/2}, \tag{14}$$

where $Z_1(s)$ is the first (say) component of $Z(s)$. Note that the LKC of $S$ depends on the local spatial correlation structure of the underlying UGRF, although we suppress this dependence in the notation $\mathcal{L}_i(S)$. The corresponding EC density in (13) is then calculated as before, but assuming the component UGRFs are isotropic. In other words, the information about the nonisotropy of the random field is transferred from the EC density to the LKC.

An explicit expression for the LKC can be found in Taylor and Adler [21], but it requires a solid grasp of differential geometry which is beyond the scope of this paper. Some idea of how the LKC can be calculated is as follows. Suppose $S$ can be embedded by a smooth transformation into a set $\tilde{S} \subset \mathbb{R}^{\tilde{N}}$ in a higher $\tilde{N} \geq N$ dimensional Euclidean space so that, in this space, local Euclidean distance is the square root of the variogram (14). Then

$$\mathcal{L}_i(S) = \mu_i(\tilde{S}).$$

The Nash Embedding Theorem guarantees that such a finite $\tilde{N}$ exists; moreover, it is bounded by $\tilde{N} \leq N(N+1)/2 + N$. However, in Taylor and Adler [21] the actual derivation of the expected EC (2) in the nonisotropic case proceeds more naturally from first principles, namely, the variogram metric (14), with isotropic random fields in Euclidean space following as a special case.

Fortunately there is a very simple expression for the $N$-dimensional LKC, which usually makes the largest contribution to the expected EC (2):

$$\mathcal{L}_N(S) = \int_S \det(\text{Var}(\dot{Z}_1(s)))^{1/2} \, ds. \tag{15}$$

Similarly, there is a simple formula for the $(N-1)$-dimensional LKC:

$$\mathcal{L}_{N-1}(S) = \tfrac{1}{2} \int_{\partial S} \det(\text{Var}(\dot{Z}_1^{\text{T}}(s)))^{1/2} \, ds,$$

where $\partial S$ is the boundary of $S$ and $\dot{Z}_1^{\text{T}}(s)$ is the derivative of $Z_1$ tangential to the boundary. The lower dimensional LKCs do not have a simple formula, except for $\mathcal{L}_0(S) = \varphi(S)$.

Of particular interest for this paper, it can be shown that in the cases we consider,

$$\mathcal{L}_i(S \times U_d) = \mu_i(S \times U_d),$$

where $S \times U_d$ is considered as a subset of $\mathbb{R}^{N+d}$. This justifies carrying on with the expected EC (2) as if $\tilde{T}(s,u)$ were jointly isotropic in $(s,u)$.

Finally, Appendix A.5 gives a neat application of the LKC version of the expected EC (13) to derive the $\chi^2$ scale space EC densities from the Gaussian scale space EC densities.

2.8. *The tube method, a direct approach to the P-value of the maximum.* The tube method, developed by Knowles and Siegmund [10], Johansen and Johnstone [9], Sun [15] and Sun and Loader [16], is a direct approach to the $P$-value of the maximum of a random field.

The approach starts by writing the random field as a Karhunen–Loéve expansion, then finding the $P$-value of the maximum of the first $m$ terms, then letting $m \to \infty$. The first $m$ terms can be written as a linear combination of $m$ spatial basis functions with $m$ independent $N(0,1)$ random variables. Conditioning on the sum of squares of these Gaussian random variables, the $P$-value of the maximum of the first $m$ terms comes down to the probability content of a subset of the $m$-dimensional unit sphere. This subset turns out to be a tube with a radius that depends on the threshold of the random field. The $P$-value calculation is thus reduced to a geometrical calculation: the ratio of the measure of the tube to the measure of the sphere. Taking expectations over the sum of squares of the Gaussian random variables completes the calculation. The final step is to let $m \to \infty$, giving a series expansion for the $P$-value similar to (2).

Sun [15] gives a general two-term series expansion for the $P$-value of the maximum of a general Gaussian random field that does not have to be homogeneous nor have a finite Karhunen–Loéve expansion, and an upper bound for an $(N-1)$-term expansion. In an unpublished manuscript by Sun (referenced in Adler [2]), it was generalized to include more general boundary cases. Sun, Loader and McCormick [17] developed simultaneous confidence regions for response curves in generalized linear models, using the tube formula, and generalized inverse Edgeworth expansions, which they developed using the Skorohod construction.

As we can see from this method, there is no need to assume isotropy, but the tube method cannot be applied to non-Gaussian fields such as $t$, $F$ or the multivariate random fields that interest us here. Moreover, in the case of Gaussian random fields, Takemura and Kuriki [19] show that the tube method's series expansion agrees with the expected EC (2) out to the number of terms in the expected EC. Thus, when they overlap, the two methods give essentially the same result, but the EC method appears to be easier to work with, and extendable to the multivariate random fields that concern us here.

**3. Random fields of multivariate test statistics.** It is not hard to see how we can apply this same union-intersection principle to other types of random fields. In Section 3.1 we shall get the EC densities for a Hotelling's $T^2$ field from the $t$ field; in Section 3.2 we shall get a $P$-value approximation (but not quite the EC density) for a Roy's maximum root field from the $F$ field; in Section 3.3 we shall get similar results for a maximum canonical correlation field from the correlation field. In fact, all the results so far known can be obtained from the correlation field, which makes the computer implementation particularly simple. This has been done in the stat_threshold function of the FMRISTAT package for the statistical analysis of fMRI data (http://www.math.mcgill.ca/keith/fmristat).

Further generalizations in Section 3.4 to multilinear forms are obvious. In fact, there is a strong link between this paper and Kuriki and Takemura [11]. The latter paper is concerned with obtaining $P$-values for maxima of multilinear forms, linear combinations of a multidimensional array (rather than just a matrix) of Gaussian random variables. Their interest is only in the $P$-value of the multilinear form itself, not a random field of multilinear forms, and their method is based on the tube method described in Section 2.8.

The same idea can easily be extended to random fields of $\bar{\chi}^2$ statistics (Lin and Lindsay [12] and Takemura and Kuriki [18]) by simply replacing the sphere $U_d$ with a cone (a subset of the sphere). This will be developed in a forthcoming paper [23]. Finally, Appendix A.5 shows how to get the $\chi^2$ scale space field from the Gaussian scale space field.

3.1. *Hotelling's $T^2$*. Hotelling's $T^2$ field is defined as

$$T(s) = \nu Z(s)'W(s)^{-1}Z(s),$$

where $W(s)$ is an independent $d \times d$ Wishart random field with $\nu$ degrees of freedom, generated as the sum of squares matrix of $\nu$ independent copies of $Z(s)$. We are interested in finding $\rho_i^{\mathrm{H}}$, the $i$-dimensional EC density of the Hotelling's $T^2$ field. Using Roy's union-intersection principle, we write the Hotelling's $T^2$ field as the maximum of the square of a (Student's) $t$ field:

$$T(s) = \max_u \tilde{T}(s,u) \qquad \text{where } \tilde{T}(s,u) = \frac{(Z(s)'u)^2}{u'W(s)u/\nu}.$$

The variance of the derivative with respect to $u$ of $Z'u$ is the identity matrix, so that $Z'u$ is a UGRF as a function of $u$ (when restricted to the unit sphere $U_d$) as well as $s$. Hence, $\tilde{T}$ is the square of a (unit) $t$ field with $\nu$ degrees of freedom.

There is a direct link between the EC of $T$ and the EC of $\tilde{T}$: for $t > 0$,

$$\varphi\{s \in S : T(s) \geq t\} = \tfrac{1}{2}\varphi\{(s,u) \in S \times U_d : \tilde{T}(s,u) \geq t\},$$

since for each fixed $s$, the excursion set $\{u \in U_d : \tilde{T}(s,u) \geq t\}$ is either empty or a pair of ellipsoidal caps with EC equal to 2. To see this, replace $u$ by $v = W^{1/2}u$, and note that the excursion set is the intersection of an ellipsoid with a pair of symmetric half-spaces. In other words, $\{(s,u) \in S \times U_d : \tilde{T}(s,u) \geq t\}$ has the same topology as two copies of $\{s \in S : T(s) \geq t\}$. Figure 3 illustrates this point for the simplest nontrivial case of $N = 1$, $d = 1$ and $\nu = \infty$.

Let $\rho_k^{\mathrm{T}}$ be the $k$-dimensional EC density of a $t$ field with $\nu$ degrees of freedom. Then

$$\mathbb{E}(\varphi\{s \in S : T(s) \geq t\}) = \tfrac{1}{2}\mathbb{E}(\varphi\{(s,u) \in S \times U_d : \tilde{T}(s,u) \geq t\})$$

$$(16) \qquad = \sum_{k=0}^{N+d} \mu_k(S \times U_d) \rho_k^{\mathrm{T}}(\sqrt{t})$$

$$= \sum_{i=0}^{N} \mu_i(S) \sum_{j=0}^{d} \mu_j(U_d) \rho_{i+j}^{\mathrm{T}}(\sqrt{t}).$$

The first step in the above (16) essentially follows from the expected EC (2), though, as noted in Section 2.7, the isotropic theory cannot be directly applied. Fortunately, as mentioned in Section 2.7, the Lipschitz–Killing curvatures of $S \times U_d$ agree with the intrinsic volumes of $S \times U_d$ considered as a subset of $\mathbb{R}^{N+d}$. In what follows, we will apply this same argument without further mention.

The last step in the above (16) follows from a result of integral geometry on the intrinsic volumes of products of sets (see Appendix A.3). Note that the factor of $\tfrac{1}{2}$ disappears because the expected EC of the excursion set of a $t$ field squared is twice that of a $t$ field not squared. Equating (16) to the expected EC (2) gives the EC density of Hotelling's $T^2$ random field:

$$\rho_i^{\mathrm{H}}(t) = \sum_{j=0}^{d} \mu_j(U_d) \rho_{i+j}^{\mathrm{T}}(\sqrt{t}).$$

Appendix A.2 gives the intrinsic volume of the sphere $\mu_j(U_d)$, and Worsley [26] gives the EC density of the $t$ field. We have thus re-derived the same result as Cao and Worsley [4], but with far less trouble.

3.2. *Roy's maximum root.* Let $V(s)$ be a Wishart random field with $\eta$ degrees of freedom and component random fields independently distributed as $Z(s)$. Let $\lambda_1(s) \geq \cdots \geq \lambda_d(s)$ be the roots of the generalized eigenvalue equation $V(s)u/\eta = W(s)u\lambda_i(s)/\nu$, where $W(s)$ is a Wishart random field with $\nu$ degrees of freedom as before. Roy's maximum root random field is $T(s) = \lambda_1(s)$. But it can also be derived from the union-intersection principle:

$$T(s) = \max_{u} \tilde{T}(s,u) \qquad \text{where } \tilde{T}(s,u) = \frac{u'V(s)u/\eta}{u'W(s)u/\nu}.$$

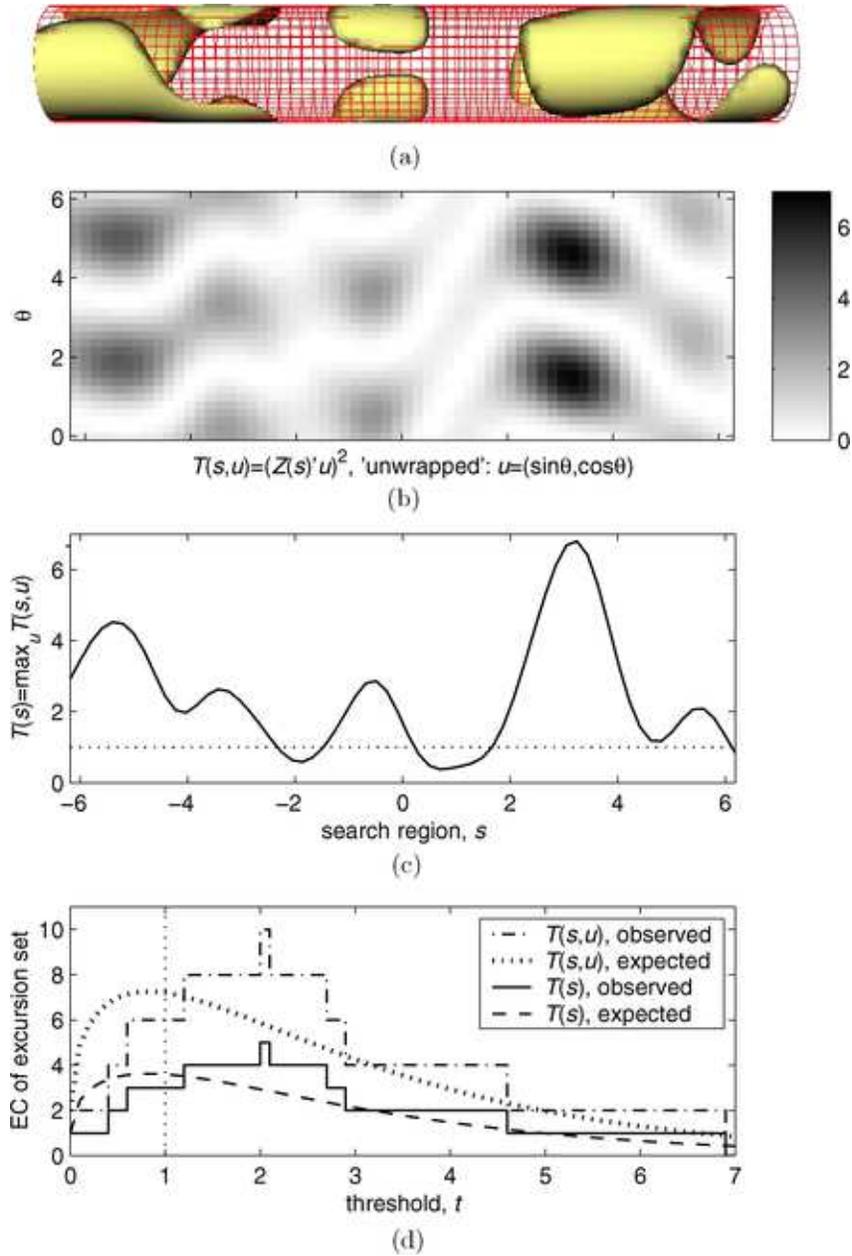

Fig. 3. *Example of Hotelling's $T^2$ field in $N=1$ dimensions, with $d=2$ components, and $\nu = \infty$ degrees of freedom. (a) Excursion set of $\tilde{T}(s,u)$ above $t=1$ (horizontal and vertical lines on graphs below); (b) "unwrapped" $\tilde{T}(s,u)$, $u=(\sin\theta,\cos\theta)$; (c) Hotelling's $T^2$, $T(s) = \max_u \tilde{T}(s,u)$; (d) observed and expected EC of $\tilde{T}$ and $T$ as a function of threshold $t$, calculated using (3) and (2). Note that the EC of $\tilde{T}$ is twice that of $T$ (e.g., at $t=1$, $6 = 2 \times 3$).*

Hotelling's $T^2$ random field is a special case with $V(s) = Z(s)Z(s)'$, that is, $\eta = 1$. Unfortunately, the EC of the excursion set of $\tilde{T}$ is no longer directly related to that of $T$, as it was in the case of Hotelling's $T^2$, as the following lemma shows:

LEMMA 1. *The EC of $\tilde{T}$ is twice the alternating sum of the EC of the roots:*

$$\varphi\{(s,u) \in S \times U_d : \tilde{T}(s,u) \geq t\} = 2\sum_{i=1}^{d}(-1)^{i-1}\varphi\{s \in S : \lambda_i(s) \geq t\}.$$

PROOF. The EC of the excursion set of $\tilde{T}$ for fixed $s$ is now more complicated, but we can find it by Morse's theorem (Morse and Cairns [13]), which states the following. Suppose $f$ is any smooth function defined on a set $A$ which takes its minimum value everywhere on the boundary of $A$, and all turning points where $\dot{f} = 0$ have nonzero $\det(\ddot{f})$ and are interior to $A$. Then Morse's theorem states that

$$\varphi(A) = \sum_{a \in A} 1_{\{\dot{f}(a)=0\}} \operatorname{sign}(\det(-\ddot{f}(a))).$$

All we have to do is choose a suitable Morse function, $f$. Let $\Lambda = \operatorname{diag}(\lambda_1, \ldots, \lambda_d)$, suppressing dependency on $s$. Then

$$\varphi\{u \in U_d : \tilde{T}(s,u) \geq t\} = \varphi\{a \in U_d : a'\Lambda a \geq t\}.$$

A suitable Morse function for $A = \{a \in U_d : a'\Lambda a \geq t\}$ is of course $f(a) = a'\Lambda a$ itself. The $2d$ turning points of $f$ in $U_d$ where $\dot{f} = 0$ are the vectors $\pm e_i$, where $e_i$ has 1 in position $i$ and 0 elsewhere. The $(d-1) \times (d-1)$ second derivative matrix $\ddot{f}(\pm e_i)$ is the diagonal matrix with elements $\lambda_j - \lambda_i$ for all $j \neq i$, and so $\operatorname{sign}(\det(-\ddot{f}(\pm e_i))) = (-1)^{i-1}$. Adding these contributions over the turning points where $f \geq t$, that is, $\lambda_i \geq t$, gives

$$\varphi\{u \in U_d : \tilde{T}(s,u) \geq t\} = 2\sum_{i=1}^{d}(-1)^{i-1}1_{\{\lambda_i \geq t\}},$$

or in other words, 2 if the number of roots greater than $t$ is odd, and 0 otherwise. Applying Morse theory to $S \times U_d$ in the same way proves the result. □

Thus in the case of Roy's maximum root, there is no direct connection between the EC of $T$ and the EC of $\tilde{T}$, except in the case of $\eta = 1$, where there is only one nonzero root (equal to Hotelling's $T^2$). For $\eta > 1$, the EC of $\tilde{T}$ is smaller than twice that of $T$, suggesting that half the EC of $\tilde{T}$ approximates the $P$-value of the maximum root maximized over $S$:

$$\mathbb{P}\left(\max_s T(s) \geq t\right) \approx \tfrac{1}{2}\mathbb{E}(\varphi\{(s,u) \in S \times U_d : \tilde{T}(s,u) \geq t\}).$$

It is difficult to give quantitative bounds on the error in the above approximation as the fields are non-Gaussian. For Gaussian fields, there exist tools such as Slepian's inequality (Adler [2]) and for the finite Karhunen–Loéve Gaussian case, the tube work of Taylor, Takemura and Adler [22]. In the non-Gaussian case, recent techniques described in Taylor, Takemura and Adler [22] give a recipe for bounding the error, though we have not applied these techniques here. If the degrees of freedom of the denominator were infinite, then the maximum of the Roy's maximum root field is a maximum of a *Gaussian* field and the results of Taylor, Takemura and Adler [22] described in Section 2.2 apply.

To evaluate the EC density, note that $\tilde{T}$ is simply an $F$ field with $\eta, \nu$ degrees of freedom and $k$-dimensional EC density $\rho_k^{\mathrm{F}}$. Then

$$\mathbb{E}(\varphi\{(s,u) \in S \times U_d : \tilde{T}(s,u) \geq t\}) = \sum_{i=0}^{N} \mu_i(S) \rho_i^{\mathrm{R}}(t),$$

where

$$\rho_i^{\mathrm{R}}(t) = \sum_{j=0}^{d} \mu_j(U_d) \rho_{i+j}^{\mathrm{F}}(t).$$

Note that $\rho_i^{\mathrm{R}}(t)$ is not the EC density of Roy's maximum root; rather, it is twice the alternating sum of the EC densities of all the roots (see Figure 4). But for high thresholds, the other roots are much smaller than the maximum, so their EC is close to zero. For this reason, we can use half $\rho_i^{\mathrm{R}}(t)$ in the approximate $P$-value (2).

It is interesting to note that the Roy's maximum root field is not always smooth. If the number of dimensions $N \geq 2$, then it can contain (with positive probability) nonsmooth local minima or "cusps" where the two largest roots are equal. Figure 5 shows an example with $N = 2$, $d = 2$, $\eta = 6$ and $\nu = \infty$. The reason is that for equality of the two roots, the two diagonal elements $v_{11}$ and $v_{22}$ of $V$ must be equal, and the off-diagonal element $v_{12}$ must be zero. These two constraints are satisfied on two zero contour lines of $v_{11} - v_{22}$ and $v_{12}$. The two lines intersect in points (with positive probability) where both constraints are satisfied, and at these points the roots are equal. This cannot happen in $N = 1$ dimensions, almost surely, so we do not see it in Figure 4. In general, equal roots, and hence cusps, will occur for any $d \geq 2$ whenever $N \geq 2$. This lack of smoothness rules out the possibility of using (8) to find the EC density of the Roy's maximum root field. This does not mean a simple expression might not exist; the conjunction random field, defined as the minimum of independent smooth random fields (Worsley and

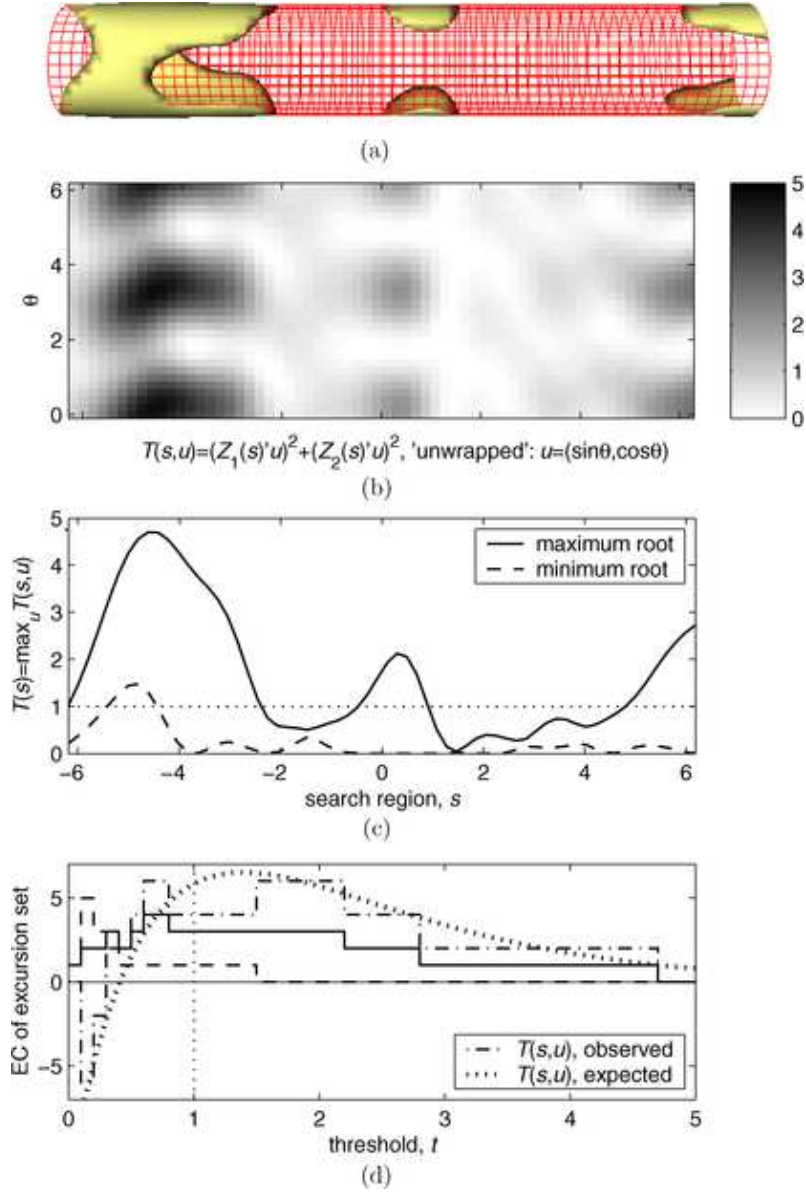

FIG. 4. *Example of Roy's maximum root field in $N = 1$ dimensions, with $d = 2$ components, and $\eta = 2$, $\nu = \infty$ degrees of freedom. (a) Excursion set of $\tilde{T}(s, u)$ above $t = 1$ (horizontal and vertical lines on graphs below); (b) "unwrapped" $\tilde{T}(s, u)$, $u = (\sin\theta, \cos\theta)$; (c) Roy's maximum root $T(s) = \max_u \tilde{T}(s, u)$ and the minimum root $\min_u \tilde{T}(s, u)$; (d) observed and expected EC of $\tilde{T}$ and observed EC of the maximum and minimum roots, as a function of threshold $t$, calculated using (3) and (2). Note that the EC of $\tilde{T}$ is twice the EC of the maximum root minus the minimum root [e.g. at $t = 1$, $4 = 2 \times (3 - 1)$]. Note that at high thresholds, the EC of the minimum root is negligible, so the EC of the maximum root is well approximated by half the EC of $\tilde{T}$.*

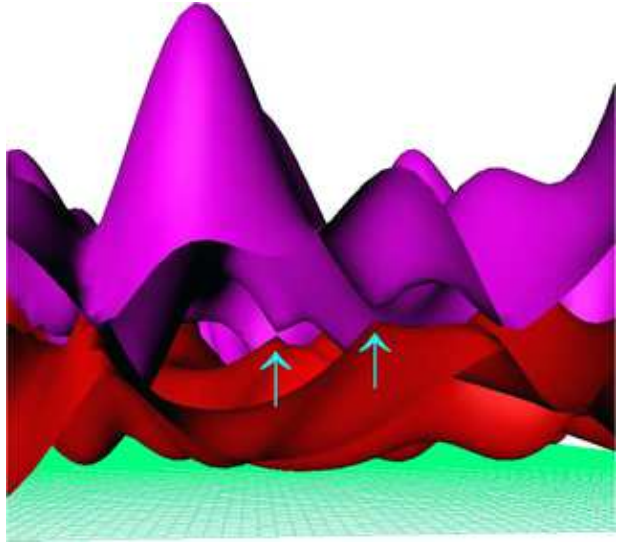

FIG. 5. *Example of cusps, or nonsmooth local minima, in a Roy's maximum root field, with $N = 2$ dimensions, $d = 2$ components, and $\eta = 6$, $\nu = \infty$ degrees of freedom. The vertical axis is the random field value (the lower lattice is zero). The upper surface is the maximum root, the lower surface is the minimum root. Two cusps are visible in the center of the picture (above arrow heads) where the upper and lower surfaces almost touch. Cusps will occur whenever $N \geq 2$ with positive probability.*

Friston [30]) also has cusps (at local maxima), yet despite this, a simple closed-form expression can be found for its EC density without using (8).

3.3. *Maximum canonical correlation.* Let $X(r)$, $r \in R \subset \mathbb{R}^M$, and $Y(s)$, $s \in S \subset \mathbb{R}^N$ be matrices of UGRFs with $c$ and $d$ columns, respectively, and the same number $n$ of rows. Let $u \in U_c$ and $v \in U_d$. Define the maximum canonical correlation random field as

$$T(r,s) = \max_{u,v} \tilde{T}(r,s,u,v)$$

$$\text{where } \tilde{T}(r,s,u,v) = \frac{u'X(r)'Y(s)v}{(u'X(r)'X(r)uv'Y(s)'Y(s)v)^{1/2}}.$$

Note that $T$ is the maximum of the canonical correlations between $X$ and $Y$, defined as the singular values of $(X'X)^{-1/2}X'Y(Y'Y)^{-1/2}$. Once again there is no direct connection between the EC of $T$ and the EC of $\tilde{T}$. Using the same approach as in Section 3.2, it can be shown that, for positive thresholds, the EC of $\{u, v : \tilde{T}(r,s,u,v) \geq t\}$ is 2 if the number of canonical correlations greater than $t$ is odd, and 0 otherwise. If $c = 1$ or $d = 1$, there is only one nonzero canonical correlation, so the EC density of $\tilde{T}$ is twice that of $T$. Otherwise, the EC of $\tilde{T}$ is smaller than twice that of $T$, suggesting that half the EC of $\tilde{T}$ approximates the $P$-value of the maximum canonical correlation maximized over $R, S$. This approximation is

$$
\begin{aligned}
&\mathbb{P}\left(\max_{r,s} T(r,s) \geq t\right) \\
&\approx \tfrac{1}{2} \sum_{i=0}^{M} \mu_i(R) \sum_{j=0}^{N} \mu_j(S) \sum_{k=0}^{c} \mu_k(U_c) \sum_{l=0}^{d} \mu_l(U_d) \rho^{\mathrm{C}}_{i+k,j+l}(t),
\end{aligned}
\tag{17}
$$

where $\rho^{\mathrm{C}}_{ij}$ is the EC density of the (cross) correlation random field $\tilde{T}$ for fixed $u, v$ given in Cao and Worsley [5]. If $M = N$, $R = S$, and the search region in $R \times S$ is confined to the "diagonal" where $r = s$, then this is called the *homologous* correlation random field (Cao and Worsley [5]), which we can generalize in the same way.

We now show that all previous results can be obtained from (17) by setting $R$ to a single point $(M = 0)$, which eliminates the need for special cases and makes it much simpler to program. For this reason, and to keep this paper self-contained, an explicit expression for $\rho^{\mathrm{C}}_{ij}$ is given in Appendix A.4. First of all, it is well known that for fixed $r, s$ the square of the maximum canonical correlation is a monotonic function of Roy's maximum root with $d$ as before and

$$
\begin{aligned}
V &= Y'X(X'X)^{-1}X'Y, & \eta &= c, \\
W &= Y'Y - V, & \nu &= n - c,
\end{aligned}
$$

and thus, Hotelling's $T^2$ ($c = 1$) and the $F$ statistic ($d = 1$) are special cases. What is not obvious is that their distributions as random fields are the same, since (17) was derived under the assumption of a random $X$ and $Y$, whereas Roy's maximum root (and the others) only require a random $Y$. However, it is easy to see that, conditional on $X$, $V$ and $W$ have the appropriate Wishart distributions with parameters that do not depend on $X$, hence, they also have the appropriate distributions marginal on $X$. Moreover, since $R$ is a point ($M = 0$), then $X$ is the same for all $s \in S$, and so the above remarks apply to the random *fields*, not just the random values at a single $s$.

3.4. *Maximum multilinear form.* Let $u_j$ be a unit $d_j$-vector, $j = 1, \ldots, D$, and $u = u_1 \otimes \cdots \otimes u_D$ be their Kronecker product, that is, a vector of all components of the form $c_1 c_2 \cdots c_D$ where $c_j$ is a component of $u_j$. Let $\tilde{Z}(s)$ be a vector of i.i.d. UGRFs of length equal to that of $u$. The maximum multilinear form random field can be defined as

$$
T(s) = \max_{u_1, \ldots, u_D} \tilde{T}(s, u_1, \ldots, u_D) \qquad \text{where } \tilde{T}(s, u_1, \ldots, u_D) = \tilde{Z}(s)' u.
$$

Note that if $D = 1$, then $T(s)^2$ is a $\chi^2$ field with $d_1$ degrees of freedom; if $D = 2$, then $T(s)^2$ is the maximum root of a $d_1 \times d_1$ Wishart random field with $d_2$ degrees of freedom, the same as $\nu$ times Roy's maximum root as $\nu \to \infty$. The EC density is straightforward:

$$\rho_N^{\mathrm{M}}(t) = \left(\prod_{j=1}^{D} \sum_{k_j=0}^{d_j} \mu_{k_j}(U_{d_j})\right) \rho_{N+k}^{\mathrm{G}}(t) \qquad \text{where } k = \sum_{j=1}^{D} k_j.$$

Noting that $T$ is unchanged if any pair of $u_j$'s are multiplied by $-1$ suggests that

$$\mathbb{P}\left(\max_s T(s) \geq t\right) \approx 2^{-(D-1)} \sum_{i=0}^{N} \mu_i(S) \rho_i^{\mathrm{M}}(t).$$

In the nonrandom field case, $N = 0$, this appears to be the same as that given by Kuriki and Takemura [11]; numerical comparisons of the coefficients of powers of $t$ are identical in every case considered so far.

## 4. Application.

4.1. *Estimating Lipschitz–Killing curvature.* Since real data is usually nonisotropic, the first step is to estimate the LKC. A method for doing this was developed in Worsley, Andermann, Koulis, MacDonald and Evans [29] and Taylor and Worsley [24], but, for completeness, we briefly describe it here as it pertains to the multivariate linear model (1).

Let $Y(s) = (Y_1(s), \ldots, Y_n(s))'$ be the $n \times d$ matrix of all the deformations data at a point $s$, and let $X = (x_1, \ldots, x_n)'$ be the $n \times p$ design matrix of the model (1), assumed of full rank. The least-squares residuals are

$$R(s) = Y(s) - X(X'X)^{-1}X'Y(s).$$

Let $R_j(s)$ be the $j$th column of $R(s)$, $j = 1, \ldots, d$. The corresponding normalized residuals are

$$Q_j(s) = R_j(s)/\|R_j(s)\|.$$

Let $e_k$ be the $N$-vector of zeros with $k$th component equal to the lattice step size along axis $k$, $k = 1, \ldots, N$. Let

$$D_j(s) = (Q_j(s+e_1) - Q_j(s), \ldots, Q_j(s+e_N) - Q_j(s))$$

be the $n \times N$ matrix of differences of $Q_j(s)$ in all lattice directions. Then $D_j(s)'D_j(s)$ is an estimator of $\mathrm{Var}(\dot{Z}_1(s))\Delta$, where $\Delta$ is the product of the lattice step sizes. Summing over all such lattice points in $S$ and averaging over components $j$, our estimator of the $N$-dimensional LKC of $S$ is

$$(18) \qquad \widehat{\mathcal{L}}_N(S) = \sum_{j=1}^{d} \sum_{s \in S} \det(D_j(s)'D_j(s))^{1/2}/d.$$

RANDOM FIELDS OF MULTIVARIATE TESTS 21
Worsley et al. [29] show that this is consistent and unbiased in the limit as the mesh size approaches zero. Note that this estimator is not invariant to a linear transformation of the components of $Y_i(s)$, and indeed, $R_j(s)$ could be replaced by any linear combination of the columns of $R(s)$ and the LKC estimator (18) would still be unbiased in the limit.

Estimating the lower dimensional LKCs is more complicated. Taylor and Worsley [24] describe a method that involves filling $S$ with a simplicial mesh (a tetrahedral mesh in $N = 3$ dimensions), a nontrivial problem. The coordinates $s \in \mathbb{R}^N$ of the mesh are then replaced by the normalized residuals $Q_j(s) \in \mathbb{R}^n$. The intrinsic volumes of all the components of the simplicial complex (points, edges, triangles, tetrahedra,...) are calculated then added together in an inclusion–exclusion type formula to give estimates of the LKCs of the union of the simplices, and hence, of $S$ itself. Again, this estimator is consistent and unbiased as the mesh size approaches zero.

In practice, this estimator produces $P$-value approximations (13) that are a little different from simply assuming that $S$ is an $N$-dimensional ball of volume $\mathcal{L}_N(S)$. The reason is that it is usually the $N$-dimensional term that makes the largest contribution to the $P$-value approximation (2). Moreover, we have already estimated $\mathcal{L}_N(S)$ quite easily by (18) without filling $S$ with a simplicial mesh. In $N = 3$ dimensions, the radius of this ball is $r = (3\widehat{\mathcal{L}}_3(S)/(4\pi))^{1/3}$ so that, from (4), the lower-dimensional LKC's can be estimated by

$$\widehat{\mathcal{L}}_0(S) = 1, \qquad \widehat{\mathcal{L}}_1(S) = 4r, \qquad \widehat{\mathcal{L}}_2(S) = 2\pi r^2.$$

This short-cut usually results in a slightly liberal $P$-value approximation since a ball has the lowest surface area (and other lower dimensional intrinsic volumes) for a given volume.

4.2. *The nonmissile trauma data.* As an illustration of the methods, we apply the $P$-value approximations for Roy's maximum root and maximum canonical correlation to the data on nonmissile trauma subjects (Tomaiuolo et al. [25]) that was analyzed in a similar way in Worsley, Taylor, Tomaiuolo and Lerch [31]. The subjects were 17 patients with nonmissile brain trauma who were in a coma for 3–14 days. MRI images were taken after the trauma, and the multivariate data were the $d = 3$ component vector deformations needed to warp the MRI images to an atlas standard (Collins et al. [6]) sampled on a 2 mm voxel lattice. The same data were also collected on a group of 19 age and sex matched controls, to give a sample size of $n = 36$. The $p = 2$ regressors were binary indicators for these two groups of subjects. Although a comparison of each trauma case with the control group might be useful in clinical applications, we follow Tomaiuolo et al. [25] and pool the trauma cases together and compare them with the control group (Tomaiuolo et al. [25] did a similar two-group comparison of white matter density, a univariate random field).

Damage is expected in white matter areas, so the search region $S$ was defined as the voxels where smoothed average control subject white matter density exceeded 5% (see Figure 1). The LKCs were estimated as in Section 4.1 by those of a ball with volume $\widehat{\mathcal{L}}_3(S) = 2571$ (radius $r = 8.5$). This choice of search region is somewhat arbitrary, but a similar search region was used in [25] for detecting group changes in white matter density. If it was felt that damage was restricted to a smaller region of higher white matter density, then the LKCs would be smaller, resulting in lower $P$-values, lower test statistic thresholds and greater power at detecting changes. However, this must be offset against the possibility that the smaller search region may have excluded regions where change really took place.

4.2.1. *Hotelling's $T^2$.* The first analysis was to look for brain damage by comparing the deformations of the 17 trauma patients with the 19 controls. We are looking at a single contrast, the difference between trauma and controls, so $\eta = 1$ and the residual degrees of freedom is $\nu = 34$. In this case Roy's maximum root is Hotelling's $T^2$. The $P = 0.05$ threshold, found by equating (13) to 0.05 and solving for $t$, was $t = 54.0$. The thresholded data, together with the estimated contrast (mean trauma—control deformations) is shown in Figure 1(a). A large region near the corpus callosum seems to be damaged. The nature of the damage, judged by the direction of the arrows, is away from the center [see Figure 1(b)]. This can be interpreted as expansion of the ventricles, or more likely, atrophy of the surrounding white matter, which causes the ventricle/white matter boundary to move outward.

4.2.2. *Roy's maximum root.* We might also be interested in functional anatomical connectivity: are there any regions of the brain whose shape (as measured by the deformations) is correlated with shape at a reference point? In other words, we add three extra regressors to the linear model whose values are the deformations at a pre-selected reference point, so now $p = 5$ and $\nu = 31$. The test statistic is now the maximum canonical correlation, or equivalently, the Roy's maximum root for these $\eta = 3$ extra regressors. We chose as the reference the point with maximum Hotelling's $T^2$ for damage, marked by axis lines in Figure 1. Figure 1(c) shows the resulting maximum canonical correlation field above the $P = 0.05$ threshold of 0.746 for the combined trauma and control data sets removing separate means for both groups. Obviously there is strong correlation with points near the reference, due to the smoothness of the data. The main feature is the strong correlation with contralateral points, indicating that brain anatomy tends to be symmetric.

A more interesting question is whether the correlations observed in the control subjects are modified by the trauma (Friston et al. [7]). In other words, is there any evidence for an interaction between group and reference vector deformations? To do this, we simply add another three regressors to the linear model whose values are the reference vector deformations for the trauma patients, and the negative of the reference vector deformations for the control subjects, to give $p = 8$ and $\nu = 28$. The resulting Roy's maximum root field for testing for these $\eta = 3$ extra regressors, thresholded at 30.3 ($P = 0.05$) is shown in Figure 1(d). Apart from changes in the neighborhood of the reference point, there is some evidence of a change in correlation at a location in the contralateral side, slightly anterior. Looking at the maximum canonical correlations in the two groups separately, we find that correlation has increased at this location from 0.719 to 0.927, perhaps indicating that the damage is strongly bilateral.

4.2.3. *Maximum canonical correlation.* If we chose to search over all reference points as well as all target points, this would lead to a maximum canonical correlation field with $X = Y$. Note that since the correlation between $X(r)$ and $Y(s)$ is the same as that between $X(s)$ and $Y(r)$ (since $X = Y$), then the $P$-value should be halved. Note also that the reference and target points must be sufficiently well separated to avoid detecting high correlation due to spatial smoothness. The parameters in our case are $M = N = c = d = 3$, and $n$ is effectively $36 - 2 = 34$ after removing the separate group means. The threshold for the maximum correlation random field at $P = 0.05$ is 0.962 from (17) and Section A.4. Computing all correlations is obviously very expensive, but aside from this, the correlation threshold of 0.962 is so high that the search over all possible correlations is unlikely to reveal much beyond the obvious symmetry reported above.

4.2.4. *Conclusion.* In conclusion, our analysis shows that damage appears to be confined to central brain regions around the ventricles, and not, as one might expect, to regions near the brain surface where the brain might have impacted the skull. Similar conclusions were reported by Tomaiuolo et al. [25] in an analysis of white matter density.

## APPENDIX

**A.1. Intrinsic volumes.** The intrinsic volume $\mu_i(A)$ of a set $A \subset \mathbb{R}^d$ with nonempty interior and bounded by a smooth hypersurface is defined as follows (Worsley [27]): $\mu_d(A)$ is the Lebesgue measure of $A$, and for $j = 0, \ldots, d-1$,

$$\mu_j(A) = \int_{\partial A} \mathrm{detr}_{d-1-j}(C)/a_{d-j},$$

where $C$ is the curvature matrix of $\partial A$, the boundary of $A$ [the $(d-1) \times (d-1)$ second derivative matrix of the inside distance between $\partial A$ and its tangent hyper-plane], $\mathrm{detr}_k(C)$ is the sum of the determinants of all $k \times k$ principal minors of $C$, and $a_k = 2\pi^{k/2}/\Gamma(k/2)$ is the Lebesgue measure of the unit $(k-1)$-sphere in $\mathbb{R}^k$. Note that $\mu_0(A) = \varphi(A)$ by the Gauss–Bonnet theorem.

**A.2. Unit sphere.** For $A = U_d$, the unit sphere in $\mathbb{R}^d$, $\mu_d(U_d) = 0$. Now $C = I_{(d-1)\times(d-1)}$ on the outside of $U_d$ and $C = -I_{(d-1)\times(d-1)}$ on the inside. Since $\mathrm{detr}_{d-1-j}(I_{(d-1)\times(d-1)}) = 2\binom{d-1}{j}$, then

$$\mu_j(U_d) = 2 \binom{d-1}{j} \frac{a_d}{a_{d-j}}$$
$$= \frac{2^{j+1}\pi^{j/2}\Gamma((d+1)/2)}{j!\Gamma((d+1-j)/2)}$$

if $d-1-j$ is even, and zero otherwise, $j = 0, \ldots, d-1$.

**A.3. Products.** We now show that

$$\mu_k(A \times B) = \sum_{i=0}^{k} \mu_i(A)\mu_{k-i}(B).$$

A very useful result of Hadwiger [8] states that if $\varphi(S)$ is a set functional that is invariant under rotations and translations, and has the additivity property

(19) $$\varphi(S_1 \cup S_2) = \varphi(S_1) + \varphi(S_2) - \varphi(S_1 \cap S_2),$$

then $\varphi(S)$ must be a linear combination of intrinsic volumes of $S$. Fixing $A$, we can see that $\varphi(B) = \mu_k(A \times B)$ has these properties. Fixing $B$ and repeating the exercise, we conclude that

$$\mu_k(A \times B) = \sum_{i=0}^{k}\sum_{j=0}^{k} c_{ij}\mu_i(A)\mu_j(B)$$

for some constants $c_{ij}$. We now determine the constants by judicious choice of $A$ and $B$. First, increasing the size of $A$ by a factor $\gamma$ increases its $k$-dimensional intrinsic volume by $\gamma^k$: $\mu_k(\gamma A) = \gamma^k \mu_k(A)$. Replacing $A, B$ by $\gamma A, \gamma B$ in the above, and noting that $\gamma A \times \gamma B = \gamma(A \times B)$, we conclude that the only nonzero constants occur when $i + j = k$. Next, let $A \subset \mathbb{R}^i$, $B \subset \mathbb{R}^{k-i}$, both with nonzero Lebesgue measure. Then $A \times B \subset \mathbb{R}^k$ has Lebesgue measure $\mu_i(A)\mu_{k-i}(B)$. Note that $\mu_j(A) = 0$ for $j > i$ and $\mu_j(B) = 0$ for $j > k - i$, so the right-hand side reduces to $c_{i,k-i}\mu_i(A)\mu_{k-i}(B)$. Since Lebesgue measure and intrinsic volume coincide in these cases, we conclude that $c_{i,k-i} = 1$. This completes the proof.

**A.4. EC density of the correlation random field.** For completeness, we reproduce the EC density $\rho_{ij}^{C}(t)$ of the correlation random field given in Cao and Worsley [5]. $\rho_{00}^{C}(t)$ is just the upper tail probability of the Beta distribution with parameters $1/2, (n-1)/2$. Let $h = i+j$. For $i>0$, $j \geq 0$ and $n > h$,

$$\rho_{ij}^{C}(t) = \frac{2^{n-2-h}(i-1)!j!}{\pi^{h/2+1}}$$

$$\times \sum_{k=0}^{\lfloor (h-1)/2 \rfloor} (-1)^k t^{h-1-2k}(1-t^2)^{(n-1-h)/2+k}$$

$$\times \sum_{l=0}^{k}\sum_{m=0}^{k} \left( \Gamma\left(\frac{n-i}{2}+l\right)\Gamma\left(\frac{n-j}{2}+m\right) \right)$$

$$\times (l!m!(k-l-m)!(n-1-h+l+m+k)!$$

$$\times (i-1-k-l+m)!(j-k-m+l)!)^{-1},$$

where $\lfloor \cdot \rfloor$ rounds down to the nearest integer, terms with negative factorials are ignored and $\rho_{ij}^{C}(t) = \rho_{ji}^{C}(t)$. The summations have been arranged for easy numerical evaluation.

**A.5. Scale space.** The Gaussian scale space random field is obtained by smoothing white noise with an isotropic spatial filter over a range of filter widths or scales, and adding the scale to the location parameters of the field (Siegmund and Worsley [14]). In essence, it is a continuous wavelet transform that is designed to be powerful at detecting localized signal of an unknown spatial scale as well as location. Let $dB(s)$ be Gaussian noise on $\mathbb{R}^N$ based on Lebesgue measure and let $f(s)$ be a filter, normalized so that $\int f^2 = 1$, and scaled so that $\int \dot{f}\dot{f}' = I_{N \times N}$. The Gaussian scale space random field with filter $f$ is defined as

$$(20) \qquad T(s,w) = w^{-N/2} \int_{\mathbb{R}^N} f((s-t)/w)\, dB(t).$$

Note that $T(s,w) \sim N(0,1)$ and $\text{Var}(\partial T/\partial s) = w^{-2N} I_{N \times N}$ at each point $s, w$. Siegmund and Worsley [14] and Worsley [28] show that for searching over a range of scales $w \in [w_1, w_2]$,

$$(21) \qquad \mathbb{E}(\varphi\{s,w \in S \times [w_1,w_2] : T(s,w) \geq t\}) = \sum_{i=0}^{N} \mu_i(S)\rho_i^{S}(t),$$

where the Gaussian scale space EC density is

$$\rho_i^{S}(t) = \frac{w_1^{-i} + w_2^{-i}}{2}\rho_i^{G}(t)$$

$$(22) \quad + \frac{w_1^{-i} - w_2^{-i}}{i} \sum_{j=0}^{\lfloor i/2 \rfloor} \frac{\kappa^{(1-2j)/2} \, (-1)^j \, i!}{(1-2j)(4\pi)^j j!(i-2j)!} \rho_{i+1-2j}^{\mathrm{G}}(t)$$

[we define $w^i/i$ as $\log(w)$ when $i = 0$]. The parameter $\kappa$ measures the variance of the derivative in the log scale direction:

$$\kappa = \int_{\mathbb{R}^N} (s'\dot{f} + (N/2)f)^2 \, ds,$$

with $\kappa = N/2$ for the Gaussian filter and $\kappa = (N+4)/2$ for the Marr filter.

The scale space random field $T(s, w)$ is nonisotropic in $(s, w)$, so the scale space result can be set in terms of the Lipschitz–Killing curvature of $S \times [w_1, w_2]$, as in (13):

$$(23) \quad \mathbb{E}(\varphi\{s, w \in S \times [w_1, w_2] : T(s, w) \geq t\}) = \sum_{i=0}^{N+1} \mathcal{L}_i(S \times [w_1, w_2]) \rho_i^{\mathrm{G}}(t).$$

Equating the two expressions (21) and (23) for the expected EC implies that $\mathcal{L}_0(S \times [w_1, w_2]) = \mu_0(S)$ and for $i \geq 1$,

$$
\begin{aligned}
&\mathcal{L}_i(S \times [w_1, w_2]) \\
&= \frac{w_1^{-i} + w_2^{-i}}{2} \mu_i(S) \\
&\quad + \sum_{j=0}^{\lfloor (N-i+1)/2 \rfloor} \frac{w_1^{-i-2j+1} - w_2^{-i-2j+1}}{i + 2j - 1} \\
&\quad \times \frac{\kappa^{(1-2j)/2}(-1)^j (i + 2j - 1)!}{(1 - 2j)(4\pi)^j j!(i - 1)!} \mu_{i+2j-1}(S).
\end{aligned}
$$
(24)

While the above expression (23) is no more compact than (21), it allows us to immediately get an expression for the EC density of the $\chi_d^2$ scale space random field with $d$ degrees of freedom, defined as the sum of squares of $d$ independent copies of the Gaussian scale space field (20). The maximum multilinear form random field can be similarly defined. Having expressed (21) in terms of Lipschitz–Killing curvatures in (23), it immediately follows that the EC density of the $\chi_d^2$ scale space field is the same as that of the Gaussian scale space field (22), but replacing Gaussian EC densities $\rho_i^{\mathrm{G}}(t)$ with $\chi_d^2$ EC densities. We have thus derived the main result of Worsley [28] with very little effort. Indeed, we have much more than this result because we can compute the EC densities of *any* scale space random field for which we can compute the EC densities in the isotropic, fixed scale, case.

Department of Mathematics and Statistics  
Stanford University  
Stanford, California 94305-4065  
USA  
and  
Départment de mathématiques  
et de statistiques  
Université de Montréal  
C.P. 6128, succ. Centre-ville  
Montréal, Québec  
Canada H3C 3J7  
E-mail: [jonathan.taylor@stanford.edu](mailto:jonathan.taylor@stanford.edu)

Department of Mathematics and Statistics  
McGill University  
Montréal, Québec  
Canada H3A 2K6  
E-mail: [keith.worsley@mcgill.ca](mailto:keith.worsley@mcgill.ca)